\documentclass[12pt]{article} \textheight 8in\textwidth
6in \oddsidemargin 0in\evensidemargin 0in
\usepackage{graphicx}
\usepackage{amsfonts}
\usepackage{epsf}

\newtheorem{theorem}{Theorem}
\newtheorem{lemma}[theorem]{Lemma}
\newtheorem{corollary}[theorem]{Corollary}

\def\pont{\hspace{-6pt}{\bf.\ }}

\newcommand{\text}[1]{\quad\mbox{#1}\quad}
\def\beq{\begin{equation}}\def\eeq{\end{equation}}
\def\beqn{\begin{eqnarray}}\def\eeqn{\end{eqnarray}}
\def\pont{\hspace{-6pt}{\bf.\ }}

\def\qed{\ifhmode\unskip\nobreak\fi\quad\ifmmode\Box\else$\Box$\fi}

\begin{document}
\title{Distribution of colors in Gallai colorings}
\author{Andr\'as Gy\'arf\'as$^a$, D\"om\"ot\"or P\'alv\"olgyi$^b$, Bal\'azs Patk\'os$^a$, Matthew Wales$^c$
\medskip \\
\small $^a$ Alfr\'ed R\'enyi Institute of Mathematics, Hungarian Academy of Sciences\\
\small P.O.B. 127, Budapest H-1364, Hungary.\\
\medskip
\small \texttt{\{gyarfas,patkos\}@renyi.hu}\\
\small $^b$ MTA-ELTE Lend\"ulet Combinatorial Geometry Research Group \\
\small Institute of Mathematics, E\"otv\"os Lor\'and University (ELTE), Budapest, Hungary\\
\medskip
\small \texttt{dom@cs.elte.hu}
\medskip\\
\small $^c$Department of Pure Mathematics and Mathematical Statistics, University of Cambridge\\
\small Wilberforce Road, Cambridge CB3 0WB\\
\small \texttt{mw637@cam.ac.uk}}

\maketitle

\begin{abstract}
A Gallai coloring is an edge coloring that avoids triangles colored with three different colors. Given integers $e_1\ge e_2 \ge \dots \ge e_k$ with $\sum_{i=1}^ke_i={n \choose 2}$ for some $n$, does there exist a Gallai $k$-coloring of $K_n$ with $e_i$ edges in color $i$? In this paper, we give several sufficient conditions and one necessary condition to guarantee a positive answer to the above question. In particular, we prove the existence of a Gallai-coloring if $e_1-e_k\le 1$ and $k \le \lfloor n/2\rfloor$. We prove that for any integer $k\ge 3$ there is a (unique) integer $g(k)$ with the following property: there exists a Gallai $k$-coloring of $K_n$ with $e_i$ edges in color $i$ for every $e_1\le\dots \le e_k$ satisfying $\sum_{i=1}^ke_i={n\choose 2}$, if and only if $n\ge g(k)$. We show that $g(3)=5$, $g(4)=8$, and $2k-2\le g(k)\le 8k^2+1$ for every $k\ge 3$.
\end{abstract}

\textit{keywords}: Gallai colorings, color class sizes, balanced colorings

\section{Introduction}

Gallai colorings (a term introduced in \cite{GYS}) of complete graphs are edge colorings that do not contain triangles colored with three different colors. For general properties of Gallai colorings see \cite{Betal,BBH,BL,GYS}. Ramsey type properties in Gallai colorings have been studied for example in \cite{GYSSS,MS}.

Here we are interested in the possible distribution of the number of edges in the color classes of Gallai colorings. For Gallai $2$-colorings there is no restriction since any $2$-coloring of $K_n$ is a Gallai coloring. In fact, a recent result of Balogh and Li \cite{BL} shows that almost all Gallai colorings of $K_n$ are $2$-colorings. For more than two colors there are restrictions and it seems not easy to find a good characterization of the sequences that are realizable as color distributions. We give a necessary and several sufficient conditions.

It is well-known that in a Gallai coloring of $K_n$ there exists a monochromatic spanning tree, consequently a color with at least $n-1$ edges. This can be generalized as follows.

\begin{lemma}\pont\label{nec} For positive integers $n,\ell$ the following is true. In every Gallai coloring of $K_n$ there exists at most $\ell$ colors such that at least $n-1+n-2+\dots+n-\ell$ edges are colored with these colors.
\end{lemma}

The lower bound of Lemma \ref{nec} is tight, we will define the construction showing its tightness in the next subsection. Lemma \ref{nec} gives a necessary condition for the distribution of colors in a Gallai coloring of $K_n$.
\begin{corollary}\pont\label{kov} Assume that $e_1\ge e_2\ge \dots \ge e_k$ and $\sum_{i=1}^k e_i={n\choose 2}$. If $K_n$ has a Gallai $k$-coloring with $e_i$ edges of color $i$, then $$\sum_{i=1}^\ell e_i\ge n-1+n-2+\dots+n-\ell$$
for every $1\le \ell \le k$.
\end{corollary}

The condition of Corollary \ref{kov} is not sufficient. For example, for $n=6$ the sequence $e_1=7,e_2=3,e_3=e_4=2,e_5=1$ satisfies the condition, but there is no Gallai coloring of $K_6$ with the given color distribution. A coloring of $K_n$ is called {\em balanced} if $e_1-e_k\le 1$. Lemma \ref{nec} (with $\ell=1$) implies that for a balanced Gallai $k$-coloring we have $k\le \lceil n/2 \rceil$.  We prove that this condition is also sufficient.

\begin{theorem}\pont\label{balanced} Assume that $k\le \lceil n/2\rceil$. Then $K_n$ has a balanced Gallai $k$-coloring.
\end{theorem}

We shall derive Theorem \ref{balanced} from the following related result.

\begin{theorem}\pont\label{div} Assume that $p\ge n-1$ and ${n\choose 2}=k\times p+q$ with $k,q\ge 0$. Then $K_n$ has a Gallai $(k+1)$-coloring with $p$ edges in $k$ color classes and $q$ edges in one color class.
\end{theorem}

Note that these theorems are equivalent if ${n\choose 2}=kp$.

\bigskip

We prove that for any fixed $k$ and large enough $n$ {\em any distribution} is possible for a Gallai $k$-coloring. Formally, there exists the function $g(k)$, the smallest integer $m$ such that for any $n\ge m$ and for any $k$-partition $e_1,\dots,e_k$ of ${n \choose 2}$ there is a Gallai $k$-coloring of $K_n$ with $e_i$ edges in color $i$. It is not quite obvious that $g(k)$ is well defined, this follows from the monotone property of Gallai colorings, stated as Lemma \ref{monotone} in Section \ref{gk}. We have the following result on $g(k)$.

\begin{theorem}\pont\label{general} For any $k\ge 3$ we have $2k-2\le g(k)\le 8k^2+1$. Moreover, $g(3)=5$ and $g(4)=8$ hold.
\end{theorem}

We pose as an open problem to determine the exact order of magnitude of $g(k)$, which we expect to be just slightly superlinear.





\subsection{General and special Gallai colorings}

A  fundamental property of Gallai colorings is the following.

\begin{theorem}\pont[Theorem A, \cite{GYS}]\label{gys}
For any Gallai coloring of $K_n$ there exist at most two colors, say $1,2$, and a decomposition of $K_n$ into $m\ge 2$ vertex disjoint complete graphs $K_{n_i}$ ($1\le i \le m$) so that all edges between $V(K_{n_i})$ and $V(K_{n_j})$ are colored with the same color and that color is either $1$ or $2$.
\end{theorem}

Briefly stated, every Gallai coloring can be obtained by a sequence of substitutions (of Gallai colored complete graphs) into a (non-trivial) $2$-colored complete graph. In our proofs we use a subfamily of Gallai colorings that we call {\em special Gallai colorings} defined as follows.

Let the vertex set of $K_n$ be $V=[0,n-1]$ and let $S(i)$ be the star with center vertex $i\in [1,n-1]$ and with edge set $\{(i,j): i>j\ge 0\}$ (note that $S(i)$ has $i$ edges). Let $T_1,\dots,T_k$ be a partition of $V\setminus\{0\}$ into nonempty parts. Then, for $j=1,2,\dots,k$, color class $C_j$ is defined as $\cup_{i\in T_j} S(i)$.  Observe that this coloring is a Gallai coloring, because for every triangle $abc$ with $n-1\ge a>b>c\ge 0$, the edges $ab,ac$ have the same color. Notice that the special coloring defined by the partition $\{n-1\},\{n-2\},\dots,\{1\}$ shows that the lower bound of Lemma \ref{nec} is sharp.

As another example, the sequence $8,3,3,1$ can be realized as a color distribution of a Gallai coloring on $K_6$ but not as a special Gallai coloring.
What we are going to prove is that the Gallai colorings claimed in Theorems \ref{balanced} and \ref{div} can be special Gallai colorings.

\section{Proof of Theorems \ref{balanced}, \ref{div}}

\noindent {\bf Proof of Lemma \ref{nec}.}
 We proceed by induction on $n$. For any $n$ and $\ell=1$ the lemma follows from the result cited above: Gallai colorings of $K_n$ contain a monochromatic spanning tree. The case $n=1$ is trivial (for any $\ell$). Let $G$ be a Gallai colored $K_n$ and $\ell\ge 2$.  Let $m$ and $e_1,e_2,\dots,e_m$ be the numbers obtained from the partition ensured by Theorem \ref{gys}. Moreover define $\alpha$ as one or two depending on the number of colors used between the $V(K_{n_i})$-s in Theorem \ref{gys}.

Let $S$ be the decreasing sequence of positive integers obtained by concatenating the sequences $S_i=e_i-1,e_i-2,\dots,1$ for $i=1,\dots,m$. For example for $$n=14, m=4, e_1=5,e_2=4,e_3=4,e_4=1$$ we get
$$S_1=4,3,2,1; S_2=3,2,1; S_3=3,2,1;S_4=\emptyset$$ and
$$S=4,3,3,3,2,2,2,1,1,1.$$ Let $S^*$ be the subsequence of $S$ defined by the first $\ell-\alpha$ elements of $S$ (if $\ell-\alpha>|S|$, then $S^*=S$). We can partition the elements of $S^*$ into sequences $S^*_i$ of length $\ell_i$, $i=1,\dots,m$ so that $\sum_{i=1}^m\ell_i=\ell-\alpha$ and the elements of $S^*_i$ form an initial segment in $S_i$ for all $i$. This is not unique, with $\alpha=2,\ell=7$ in the example above we have $$S^*=4,3,3,3,2$$ and we can have $$S_1^*=4,3; S_2^*=3; S_3^*=3,2; S_4^*=\emptyset \mbox{  or  } S_1^*=4,3,2; S_2^*=3; S_3^*=3; S_4^*=\emptyset$$.

Now we can apply induction for each Gallai colored $K_{e_i}$ with $\ell_i$ and together with the $\alpha$ colors from the color set $\{1,2\}$ we get at most $\ell$ colors so that the number of edges in these colors is at least $$t=\sum_{1\le i<j \le m} e_ie_j + \sum_{i=1}^m ||S_i^*||$$
where $||S_i^*||$ denotes the sum of elements of the sequence $S_i^*$. Since the elements of $S^*_i$ form an initial segment in $S_i$ for all $i$, the number of edges of $G$ not among the $t$ selected edges is at most the sum $s$ of the $(n-m)-(\ell-\alpha)$ elements of $S^*$ that are not in the $S_i^*$-s. But $$(n-m)-(\ell-\alpha)\le n-\ell-1$$ because if $m=2$, then $\alpha=1$ and if $m\ge 3$, then $\alpha\le 2$. Therefore $s$ is at most the sum of the first  $n-\ell-1$ numbers from the set $\{1,2,\dots ,n-1\}$ and
$$t\ge {n\choose 2}-(1+2+\dots+n-\ell-1)=n-1+n-2+\dots+n-\ell$$ proving the lemma.  \qed

\medskip

\noindent {\bf Proof of Theorem \ref{div}.}
We prove by induction on $n$ that the required Gallai coloring can be chosen to be a special Gallai coloring.
The proof is built around three primary `moves' we can make to reduce to a smaller case. When $p\leq 2n-3$ we form some classes out of a pair of stars $S(i)$ and $S(j), i\neq j$.
	If $p\geq 2n-2$, we can't form classes from a pair of stars. We will still place pairs of stars in each class of size $p$, but then apply induction with the same value of $k$, and smaller values of $p$ and $n$.
	To resolve a final case, we will also remove a star from the class of size $q$.

	We proceed by strong induction on $n$. The base case $n=1$ is trival. We can assume $q\leq p-1$ by finding a partition into one part of size $q \pmod p$, and $k + \lfloor \frac{q}{p}\rfloor$ parts of size $p$. Observe the hypotheses are still satisfied, and at the end we can merge these newly created parts to form a single size q part.
	If $q\geq n-1$, we can find a Gallai coloring of $K_{n-1}$ with one part of size $q-(n-1)$ and $k$ classes of size $p$. Adding a star on a new vertex in the final color class gives the desired Gallai coloring.
	We break into 2 primary cases, depending on whether removing a pair of stars can reduce the number of classes. Thus we can assume

\begin{equation}\label{qpn}
q\le n-2, q\le p-1
\end{equation}
	
	\textbf{Case A:} $p\leq 2n-3$
	
	Note that $0\leq (p-n+1)<(n-1)$. We will form some classes out of pairs of stars.
	If $p$ is odd, we form classes of size $p$ from each of the star pairs $(S(n-1),S(p-n+1)),...,(S(\frac{p+1}{2}),S(\frac{p-1}{2}))$. What remains is a complete graph on $p-n+1$ vertices. By induction, since $p$ has not changed (and $n$ decreases), we can construct a Gallai coloring of this graph into some classes of size $p$, and one class of size $q$. Adding the stars sequentially gives the desired Gallai coloring.
	
	As such, we may assume $p$ is even. We form size $p$ classes from each of the star pairs $(S(n-1),S(p-n+1)),...,(S(\frac{p}{2}+1),S(\frac{p}{2}-1))$. What remains is a complete graph on $p-n+1$ vertices, and a star at $\frac{p}{2}$, consisting of $\frac{p}{2}$ edges.
	
	We construct a Gallai Coloring of the $K_{p-n+1}$ into some parts of size $\frac{p}{2}$, and one of size $q$ by induction - recall ${n\choose 2} = kp + q$, and we only removed classes of size $p$ and a star of size $\frac{p}{2}$. Since $p\leq 2n$ (which implies $\frac{p}{2}\leq p-n$) the induction hypotheses hold. We form a final size $\frac{p}{2}$ class from the leftover star, and pairing off the size $\frac{p}{2}$ classes arbitrarily gives the desired Gallai coloring of $K_n$.
	
	This resolves Case A.
	
	\textbf{Case B:} $p\geq 2n-2$
	
	In this situation, we try and use the second move to reduce to a smaller case. Observe that $k\leq\frac{n}{4}$ because ${n\choose 2}\geq kp$. We can therefore form pairs of stars $S(n-i),S(n-2k-1+i)$, with a combined $2n-2k-1$ edges. These are all distinct stars for $1\leq i \leq k$.
	With this in mind, let $p' = p-2n+2k+1$, and $n' = n-2k$. We use these for the induction.
	
	If $p'\geq n'-1$, then by induction we can find a Gallai coloring of $K_{n'}$ with $k$ classes of size $p'$ and one class of size $q$. We then add the star pairs to the size $p'$ classes to obtain the desired Gallai coloring.
	So, assume this is not the case, and therefore $p\leq 3n-4k-3$. Since also $p\geq 2n-2$, we deduce $4k\leq n-1$. Let $4k=n-\delta$, with $\delta\in\mathbb{Z}^{\geq 1}$. Since by (\ref{qpn}) $q\leq n-2$,

$${n\choose 2}\leq k(3n-4k-3) + (n-2).$$
	
	Substituting $4k = n-\delta$, we obtain
$$	n\delta + \delta(\delta-3)\leq 3n-4.$$
	
	This implies that $\delta < 3$, and so either $\delta = 1$ or $\delta = 2$. Since $2n-2 \leq p \leq 2n+\delta - 3$, there are very few remaining possibilities. We check these individually.
	
	If $\delta = 1$, $n = 4k+1$ and $p = 2n-2 = 8k$. Since ${n\choose 2} = kp + q$, we need $q = 2k$. Remove the k pairs of stars as earlier from the size $p$ classes. It remains to find a Gallai coloring of $K_{n'} (n' = 2k+1)$ with $k$ classes of size $p' = 2k-1$ and one class of size $2k$. We remove a star in the size $q$ class, and now need to Gallai color $K_{2k}$ with $k$ classes of size $2k-1$. But we can do this by induction, and adding back stars gives the desired result. In fact, this case can be settled directly too, by taking one triple of stars $S(4k),S(4k-1),S(1)$ and $k-1$ quadruples of stars $S(4k-2i),S(4k-2i-1),S(2i),S(2i+1)$ (for $i=1,\dots,k-1$). The only unused star $S(2k)$ gives the part with $q$ edges.
	
	Thus $\delta = 2$, and so $n = 4k+2$ and either $p = 8k+2$ or $p = 8k+3$. If $p = 8k+2$, we find $q = 4k+1 = n-1$. But this contradicts that $q\leq n-2$. Therefore, $p = 8k+3$, and hence $q = 3k+1$. But in this case, we remove a pair of stars, as at the start of Case B, from each size $p$ class. It'ls now sufficient to find a coloring of $K_{n'} ,n' =  n-2k = 2k+2$, with $k$ classes of size $p' = 2k$ and 1 class of size $q = 3k+1$. We remove a star of size $n'-1$ from the size $q$ class, and we then need a Gallai coloring of $K_{2k+1}$ with $k$ classes of size $2k$ and 1 class of size $k$. But this exists by induction. Adding the stars gives the desired coloring, and completes Case B.  \qed

\noindent {\bf Proof of Theorem \ref{balanced}.}
Observe that it is enough to prove Theorem \ref{balanced}  for the case $n/4\le k$. Indeed, otherwise with a suitable positive integer $r$ we have $n/4 \le rk \le \lceil n/2\rceil$ and Theorem \ref{balanced} provides special balanced Gallai coloring with $rk$ colors. Then, grouping the colors into $k$ parts  so that every part consists of $r$ old color classes, we get the required solution. (Note that a special Gallai coloring remains special after merging some color classes.) Therefore we can assume $n=2k+i$ with $0\le i\le 2k$ and write ${i\choose 2}=\ell k +m$ where $0\le m \le k-1$.
In a balanced distribution, there are $k-m$ color classes with $Z=2k+2i+\ell-1$ edges and $m$ color classes with $Z'=2k+2i+\ell$ edges.
Assume that $\ell$ is even. We can create $k-{\ell\over 2}$ pairs of stars as follows:
\begin{equation}\label{part}
(S(2k+i-1),S(i+\ell)),\dots. (S(k+i+{\ell\over 2}),S(k+i+{\ell\over 2}-1)).
\end{equation}

If $k-m\ge k-{\ell\over 2}$, then the $k-{\ell\over 2}$ pairs of centers in (\ref{part}) cover all numbers in $[i+\ell,n-1]$ exactly once and each pair of stars has total size $Z$. This implies that
${i+\ell\choose 2}=({\ell\over 2}-m)Z+mZ'$.  Thus, by induction, we can find a balanced special Gallai coloring into ${\ell\over 2}$ parts, ${\ell\over 2}-m$ color classes of size $Z$ and $m$ color classes with size $Z'$. Thus together with the pairs in (\ref{part}) we have the required coloring on $K_n$. Since ${\ell\over 2}\le {i+\ell\over 2}$ and $i+\ell < n=2k+i$, the induction is justified.

If $k-m< k-{\ell\over 2}$, then we have too many pairs of stars with total size $Z$, thus we need to stop earlier in the pairings at (\ref{part}). We make the following modification in the pairing.

\begin{equation}\label{partmod1}
(S(2k+i-1),S(i+\ell)),\dots, (S(k+i+m),S(k+i-m+\ell-1))
\end{equation}

defining $k-m$ pairs with sum $Z$ and continue with pairings with sum $Z'$ as follows.

\begin{equation}\label{partmod2}
(S(k+i+m-1),S(k+i-m+\ell+1)),\dots, (S(k+i+{\ell\over 2}+1),S(k+i+{\ell\over 2}-1))
\end{equation}
defining $m-{\ell\over 2}-1$ pairs with sum $Z'$.  Observe that we did not use two numbers as centers of stars from the interval $[i+\ell,n-1]$, namely we skipped $k+i-m+\ell$ and at the end $k+i+{\ell\over 2}$ was not used in the pairings. The sum of these  is $Z'-(m-{\ell\over 2})$. This implies that ${i+\ell\choose 2}={\ell\over 2}Z'+m-{\ell\over 2}$. Since $Z'=2k+2i+\ell\ge i+\ell-1$, Theorem \ref{div} can be applied with $i+\ell$ in the role of $n$, ${\ell\over 2}$ in the role of $k$, $Z'$ in the role of $p$, $m-{\ell\over 2}$ in the role of $q$ to get a special Gallai coloring with ${\ell\over 2}$ parts of size $Z'$ and one part of size $m-{\ell\over 2}$. Adding the star-pairs from (\ref{partmod1}), (\ref{partmod2}) together with the two stars unused at (\ref{partmod2}) we have the required coloring.

Suppose now $\ell$ is odd. If $k-{\ell+1\over 2}\le m$, then we start with $k-{\ell+1\over 2}$ pairs of stars with sum $Z'$ with centers in $[i+l+1,n-1]$.
\begin{equation}\label{partmod3}
(S(2k+i-1),S(i+\ell+1)),\dots, (S(k+i+{\ell+1\over 2}),S(k+i+{\ell-1\over 2})).
\end{equation}

  We need $m-k+{\ell+1\over 2}$ further parts of size $Z'$ and $k-m$  parts of size $Z$ with center in $[1,i+\ell]$. This can be done by induction.

  Otherwise, when $k-{\ell+1\over 2}> m$,
we can define $m$ pairs of sum $Z'$ and $k-m-{\ell+1\over 2}-1$ pairs of sum $Z$ as follows:

\begin{equation}\label{partmod4}
(S(2k+i-1),S(i+\ell+1)),\dots, (S(2k+i-m),S(i+\ell+m)),
\end{equation}

\begin{equation}\label{partmod5}
(S(2k+i+m-2),S(i+\ell+m+1)),\dots, (S(k+i+{\ell+1\over 2}),S(k+i+{\ell-1\over 2}-1)).
\end{equation}

Note that in  (\ref{partmod4}), (\ref{partmod5}) we get stars with all centers from $[i+\ell+1,n-1]$ except the ones with centers $2k+i-m-1$ and $k+i+{\ell-1\over 2}$. The size of these together is $Z+k-m-1-{\ell+1\over 2}$
therefore we have ${i+\ell+1\choose 2}= ({\ell+1\over 2}-1)Z+m-k+1+{\ell+1\over 2}$. Since $Z=2k+2i+2\ell-1\ge i+\ell-1$, Theorem \ref{div} can be applied to get a special Gallai coloring on $K_{i+\ell}$ with ${\ell+1\over 2}-1$ parts of size $Z$  and one part of size $m-k+1+{\ell+1\over 2}$. Adding the pairs at (\ref{partmod4}), (\ref{partmod5}) plus the two exceptional stars, we get the desired balanced Gallai $k$-coloring. \qed

\section{Bounds on $g(k)$}\label{gk}

A simple general lower bound is $g(k)\ge 2k-2$, shown by the color distribution ${2k-3\choose 2}-k+1,1,\dots,1$ for $K_{2k-3}$. Indeed, this distribution is impossible for a Gallai $k$-coloring, since among the edges of color $2,\dots,k$ there must be two intersecting ones defining a multicolored triangle. This enables us to show that Gallai-colorability is a monotone property in the following sense, implying that $g(k)$ is well defined.

\begin{lemma}\pont\label{monotone} Assume that for some $n$, $K_n$ has a Gallai $k$-coloring for every distribution $e_1\ge\dots \ge e_k$ satisfying $\sum_{i=1}^k e_i={n\choose 2}$. Then this statement remains true for $n+1$ as well.
\end{lemma}

\noindent \textbf{Proof.}
Assume that $K_n$ has a Gallai $k$-coloring for every distribution with $\sum_{i=1}^k e_i={n\choose 2}$.
As shown by the distribution ${n\choose 2}-k+1,1,\dots,1$, we must have $n\ge 2k-2$, i.e. ${n+2\over 2}\ge k$.
Suppose that there is some distribution $S$, $e_1\ge \dots \ge e_k$ with $\sum_{i=1}^k e_i={n+1\choose 2}$ for which there is no Gallai $k$-coloring of $K_{n+1}$. We have $e_1\ge n$ otherwise
$$\sum_{i=1}^k e_i\le k(n-1)\le {n+2\over 2}(n-1)< {n+1\choose 2},$$
a contradiction.  Replacing $e_1$ with $e=e_1-n$ and keeping the other $e_i$-s, we have a distribution realizable as a Gallai coloring on $K_n$. Adding a star to $K_n$ with $n$ edges in the color of $e$, we get a Gallai $k$-coloring on $K_{n+1}$ with distribution $S$, a contradiction.  \qed

\bigskip

\noindent {\bf Proof of Theorem \ref{general}.}
First we prove the general upper bound $g(k)\le 8k^2+1$. We proceed by induction on $k$ with the base case given by the $k=3$ part of the theorem.

Let $n=8k^2+1$.
Because of Lemma \ref{monotone} it is enough to prove that for any $\sum_{i=1}^k e_i={n\choose 2}$ there is a Gallai $k$-coloring of $K_{n}$.
We will use $e_1\ge {n\choose 2}/k \ge 32k^3$ to eliminate color $k$ and reduce the problem to $k-1$ colors.

We give a procedure to color the edges.
Denote the current number of vertices by $n'$ (initially $n'=n$) and the current number of required edges that need color $i$ by $e_i'$, so that we will always have $\sum_{i=1}^k e_i'={n'\choose 2}$.
If $e_k'\ge n'-1$, put the star $S(n')$ in color $k$ into the graph, reduce $e'_k$ by $n'$ and reduce $n'$ by one.
We repeat this until $e_k'<n'-1$. As $e_k\le {n \choose 2}/k\le 32k^3+4k\le 5k(8k^2-5k+1)=5k(n-5k)$ holds if $k\ge 5$ (in fact if $k\ge 3$), we placed at most $5k$ stars of color $k$ in this first phase.

Since at this point $e_k'<n$, all the further needed edges of color $k$ can be placed within the next set $A$ of $4k+1$ vertices of $K_{n'}$, as $n=8k^2+1<{4k+1 \choose 2}$.
We color all other edges adjacent to $A$ with color 1. There are at most $(4k+1)\cdot 8k^2-{4k+1\choose 2}\le 32k^3\le e_1$ such edges, thus this second phase is indeed doable.

Finally we need to check if after removing the vertices of the stars and the vertices in $A$,  we are still left with at least $8(k-1)^2+1$ vertices. Or equivalently, we removed less than $16k$ vertices.  Indeed, we removed at most $5k+4k+1<16k$ vertices. Because we eliminated all the required $e_k$ edges of color $k$, by induction we can give a Gallai $(k-1)$-coloring on the remaining complete graph $K$ according to the sequence $e_1',e_2,\dots,e_{k-1}$ where $e_1'$ is the remaining number of color 1 edges to be placed. The obtained $k$-coloring is a Gallai coloring because the first two phases give a $2$-coloring and every vertex is homogeneously connected in these colors to the coloring of $K$.  \qed

\medskip

Next we prove $g(3)=5$. By Lemma \ref{monotone} it is enough to consider the case $n=5$. There are eight possible distributions, six of them with a straightforward special Gallai coloring:
\begin{equation}\label{small}
(7,2,1):S_4\cup S_3,S_2,S_1; (6,3,1):S_4\cup S_2,S_3,S_1; (5,4,1): S_3\cup S_2,S_4,S_1;
\end{equation}
\begin{equation}
(5,3,2): S_4\cup S_1,S_3,S_2; (4,3,3): S_4,S_3,S_2\cup S_1; (4,4,2): S_4,S_3\cup S_1,S_2
\end{equation}
The distribution $(8,1,1)$ is realized by taking two vertex disjoint edges in colors one and two and color all remaining edges by the third color. The distribution $(6,2,2)$  can be realized by taking a $K_{2,3}$ in color one and the other two colors take care of themselves.
\bigskip

We finish by proving $g(4)=8$. It is left to the reader to check that no Gallai $4$-coloring exists on $K_7$ with color distribution $(9,4,4,4)$, thus $g(4)\ge 8$.

By Lemma \ref{monotone}, we may assume that $n=8$. Note that induction on $k$ may work at this point: if some $e_i=7$, then we can delete $e_i$ and we get a sequence of three numbers whose sum is $21$ thus by the $k=3$ part of the theorem, there is a Gallai $3$-coloring on $K_7$ with the given distribution. Then extending this with a star in the fourth color, we get the required coloring on $K_8$. This idea can be carried further as follows. If there is no $7$ in the sequence but there is a $6$, then we can replace $e_1$ by $e_1-7$ and delete the $e_i$ with value $6$ and get three numbers whose sum is $15$ and by the $k=3$ part of the theorem there is a Gallai $3$-coloring on $K_6$ with the given distribution. Also, if there is no $7,6$ but there is a $5$ in the sequence, then we can reduce the largest or the two largest numbers by $7+6=13$ to get a sequence with sum $15$. Deleting the $e_i$ with value $5$, we can apply  the $k=3$ part of the theorem with $K_5$ and applying the corresponding extension.

Excluding the values $7,6,5$ from the sequence $e_1\ge e_2\ge e_3\ge e_4$ we have one, two or three elements of the sequence larger than $7$.

\noindent
{\bf Case 1. $e_1\ge 8>e_2$.} In this case $3\le e_2+e_3+e_4 \le 12$. We show that  the distribution $e_2,e_3,e_4$ can be realized as a Gallai coloring on the union of  vertex disjoint complete graphs.  Then the complement of this graph has $e_1$ edges and form a complete partite graph on eight vertices, providing the required Gallai $4$-coloring.

We give the partition according to the sum $S=e_2+e_3+e_4$ and show only the non-trivial part of the partitions (the $K_1$ parts are omitted). Finding the corresponding Gallai colorings is easy, we leave this to the reader.
\begin{itemize}
\item $S=12$.  $(4,4,4)\rightarrow K_4 \cup K_4$
\item $S=11$.  $(4,4,3)\rightarrow K_5 \cup K_2$
\item $S=10$.  $(4,4,2),(4,3,3)\rightarrow K_5$
\item $S=9$.  $(4,4,1),(4,3,2),(3,3,3)\rightarrow K_4 \cup K_3$
\item $S=8$.  $(4,3,1),(4,2,2),(3,3,2)\rightarrow K_4 \cup K_2\cup K_2$
\item $S=7$.  $(4,2,1),(3,3,1),(3,2,2)\rightarrow K_4 \cup K_2$
\item $S=6$.  $(4,1,1),(3,2,1),(2,2,2)\rightarrow K_3 \cup K_3$
\item $S=5$.  $(3,1,1),(2,2,1)\rightarrow K_3 \cup K_2\cup K_2$
\item $S=4$.  $(2,1,1)\rightarrow K_3 \cup K_2$
\item $S=3$.  $(1,1,1)\rightarrow K_2 \cup K_2\cup K_2$
\end{itemize}

\noindent
{\bf Case 2. $e_1,e_2\ge 8>e_3$.}  Replacing $e_2$ by $e_2-7$ and reordering, we have a new sequence with only $e_1\ge 8$.  Now the method of Case 1 can be used with $n=7$. From easy inspection of $e_2$ and $e_3+e_4\le 8$ (before the reduction), we get $S=e_2+e_3+e_4 \le 11$. Since only the case $S=12$ used eight vertices we finish as in Case 1.

\noindent
{\bf Case 3. $e_1\ge e_2 \ge e_3\ge 8>e_4$.} Replacing $e_2,e_3$ by $e_2-7,e_3-6$ and reordering, we have a new sequence with only one element larger than $7$. Now $n=6$ and easy inspection shows that $S=e_2+e_3+e_4\le 7$,
when case 1 uses more than six vertices only for $S=5$. Thus the only problem is when we end up with a sequence with $e_2+e_3+e_4=5$. However, this is impossible because then $e_1=10$ and from the assumption of the subcase $e_2\ge e_3\ge 8$ was true before the reduction, thus $e_4\le 0$, contradiction.  

\medskip




\bigskip

\bigskip

\textbf{Funding}: Research supported by the National Research, Development and Innovation Office - NKFIH under the grants SNN 129364 and K 116769, and by the Lend\"ulet program of the Hungarian Academy of Sciences (MTA), under grant number LP2017-19/2017.

\medskip

\end{document}